
\baselineskip=14pt
\parskip=10pt

\font\eightrm=cmr8 

\magnification=\magstephalf

\def\1{{\overline{1}}}
\def\2{{\overline{2}}}
\parindent=0pt
\overfullrule=0in

\def\frac#1#2{{#1 \over #2}}
\centerline
{\bf Automatic Counting of Generalized Latin Rectangles and Trapezoids}
\bigskip
\centerline
{\it George SPAHN and Doron ZEILBERGER}

\bigskip

{\bf Abstract}: In this case study in ``fully automated enumeration'', we illustrate how to take full advantage of symbolic computation by developing  (what we call) `symbolic-dynamical-programming' algorithms
for computing many terms of `hard to compute sequences', namely the number of  Latin trapezoids, generalized derangements, and generalized three-rowed Latin rectangles.
At  the end we also  sketch the proof of a generalization of
Ira Gessel's 1987 theorem that says that for any number of rows, k, the number of Latin rectangles with k rows and n columns is P-recursive in n.
Our algorithms are fully implemented in Maple, and generated quite a few terms of such sequences.

{\bf How it all Started}

Last year, for a few months, the {\it New York Times} magazine published a puzzle, created  by Wei-Hwa Huang,  called {\it Triangulum}. One is given a discrete equilateral triangle
consisting of $1+2+3+4+5=15$ empty circles, where the bottom row has 5 empty circles, the next row has four empty circles, $\dots$ and the top (fifth) row has one empty circle.
Intertwined between these circles are $15$ mostly empty triangles, but a few of them are filled with integers. The solver has to fill-in the empty circles with the
integers $1,2,3,4,5$ in such a way that each of the five horizonal lines, and each of the $10$ diagonal lines have {\bf distinct} entries, and in addition
the labels in the entries in the circles around each of the non-empty triangles add-up to the number in that triangle. See (and play!)

{\tt https://sites.math.rutgers.edu/\~{}zeilberg/EM21/projects/LTgame.html} \quad .

Of course, from a computational point of view, these puzzles are trivial, and one can do it easily by brute force. It is easy to see that there are only
$4$ such {\it reduced} configurations where the bottom row is $12345$.  Here there are:

$$
\matrix{   &   &    &    &  2  &     &      &     &     \cr
           &   &    & 3  &     &  4  &      &     &     \cr
           &   & 4  &    &  5  &     &   1  &     &     \cr
           & 5 &    & 1  &     &  2  &      &  3  &     \cr
         1 &   & 2  &    &  3  &     &  4   &     &  5 } \quad, \quad
\matrix{   &   &    &    &  3  &     &      &     &     \cr
           &   &    & 5  &     &  1  &      &     &     \cr
           &   & 2  &    &  3  &     &   4  &     &     \cr
           & 4 &    & 5  &     &  1  &      &  2  &     \cr
         1 &   & 2  &    &  3  &     &  4   &     &  5 } \quad, \quad
$$

$$
\matrix{   &   &    &    &  3  &     &      &     &     \cr
           &   &    & 2  &     &  4  &      &     &     \cr
           &   & 5  &    &  3  &     &   1  &     &     \cr
           & 4 &    & 1  &     &  5  &      &  2  &     \cr
         1 &   & 2  &    &  3  &     &  4   &     &  5 } \quad, \quad
\matrix{   &   &    &    &  4  &     &      &     &     \cr
           &   &    & 2  &     &  3  &      &     &     \cr
           &   & 5  &    &  1  &     &   2  &     &     \cr
           & 3 &    & 4  &     &  5  &      &  1  &     \cr
         1 &   & 2  &    &  3  &     &  4   &     &  5 } \quad, \quad
$$

hence altogether there are $4 \cdot 5!=480$ legal ways, and then even a very
stupid, but patient, human, who knows how to add, can pick which of the $480$ possibilities meet the extra conditions.

Let's call such equilateral discrete triangles {\it Latin Triangles}, and consider the problem of finding the {\it exact} number of
Latin triangles with side-length $n$. By brute force, one can easily get the first  few terms. Deciding that the bottom row is $1....n$, then
for $n=3,4,5,6,7$ these numbers are $1,0, 4, 236, 27820$, respectively (to get the total number , multiply by $n!$).
This sequence was not (Aug. 13, 2021) in the OEIS.

To see all the  reduced  Latin triangles of side-length up to $7$, see:

{\tt https://sites.math.rutgers.edu/\~{}zeilberg/tokhniot/oLatinTrapezoids2.txt} \quad .

We are almost sure that human-kind
will never know the exact number of Latin triangles with side-length $30$. For the apparently much easier problem of counting Latin squares,  given by OEIS sequence A2860
({\tt https://oeis.org/A002860}), currently only $11$ terms are known.

So let's settle for something {\it easier}, but far from trivial: Count Latin trapezoids with three rows. here is an example with base-length $9$:

$$
\matrix{  
           &    & 6  &    &  3  &     &   7  &     & 8  &   & 2   &   &  5   &    & 1   &  &    \cr
           &  8  &   &  1  &    &  9   &     &  2   &   & 7  &    & 3  &     &  4  &    & 6 &    \cr
         1  &   & 2  &    &  3  &     &  4   &     &  5 &  &  6  &  &  7  &  &  8  &  &  9 
} \quad .
$$

Even counting Latin trapezoids with $3$ rows is already challenging, and using the methods of this article
we were able to find the first $100$ terms. Here are the first few terms, starting with $n=3$:
$$
1, 6, 68, 1670, 67295, 3825722, 285667270, 26889145828, 3102187523467, 
429700007845870, 
$$
$$
70303573947346474, 13405343287124139802, 2945521072579394529097, 
$$
$$
738633749151050116349946, 209620243382776121032416188, \dots
$$
The remaining terms may be viewed here:

{\tt https://sites.math.rutgers.edu/\~{}zeilberg/tokhniot/oLatinTrapezoids1.txt} \quad .

This sequence was also not in the OEIS (viewed Aug. 13, 2021).

To see all the  $3$-rowed Latin trapezoids with base-lengths from  $3$ to $7$, see:

{\tt https://sites.math.rutgers.edu/\~{}zeilberg/tokhniot/oLatinTrapezoids3.txt} \quad .

All the above date was obtained using the Maple package {\tt LatinTrapezoids.txt}, available from

{\tt https://sites.math.rutgers.edu/\~{}zeilberg/tokhniot/LatinTrapezoids.txt} \quad .

After we wrote that Maple package, we realized that Latin Trapezoids may be viewed as a kind of {\it generalized Latin Rectangle}, so we went on to write a much more general Maple package, called
{\tt GenLatinRecs.txt}, available from:

{\tt https://sites.math.rutgers.edu/\~{}zeilberg/tokhniot/GenLatinRecs.txt} \quad .

The rest of this article will consist in explaining the ideas and theory behind this very versatile Maple package, that can crank out many terms for the enumerating sequences for many families
of $3$-rowed arrays with many restrictions. The classical case of counting $3$-rowed Latin rectangles, is OEIS sequence  A186, ({\tt https://oeis.org/A000186}) .
In the present article we generalize it considerably.

{\bf What is a Latin Rectangle?}

A $k \times n$ Latin rectangle is a $k \times n$ array of integers where every one of the $k$ rows is a permutation of $\{1,2, \dots, n \}$, and the entries of
each column are {\bf distinct}. In other words, the array $(M_{i,j}), (1 \leq i \leq k, 1 \leq j \leq n) $ is a {\bf Latin rectangle} if
the following conditions are satisfied

$\bullet$ Each row $i$, $1 \leq i \leq k$, $M_{i1} \dots M_{in}$ is a permutation of $\{1,2,\dots,n\}$ .

$\bullet$ For two different rows $i$ and $i'$,   $M_{ij} \neq M_{i'j}$ for all $1 \leq j \leq n$.

It is convenient to declare that the first row is the identity permutation $12 \dots n$, and such creatures are called {\it reduced} Latin rectangles. Of course in order
to find the total number, all one has to do is multiply by $n!$.

The problem of enumerating $ 2 \times n$ reduced Latin rectangles is nothing but the famous {\it probl\`eme des rencontres}, due to Montmort,
(see [R], p. 58) , and we will soon revisit this as a {\it motivating case}. Things get much harder for $k>2$. As mentioned above, the case $k=3$ is OEIS sequence A186, and
has been treated quite extensively (see the references there), in particular see [BL] [G1]. Things get {\it really} harder, computationally, for $k>3$, but a beautiful
{\it  theoretical} result of Ira Gessel [G2] asserts that for each specific (numeric) $k$, the enumerating sequence (in $n$), is $P$-recursive, or equivalently, the
generating function is $D$-finite, i.e. satisfies a linear differential equation with polynomial coefficients.

In this article we will develop an efficient algorithm to enumerate as many terms as possible for counting {\it generalized Latin rectangles}. The actual implementation
is only carried for two and three rows, but in the last section we will indicate how this can be extended for any number of rows, and show how the
algorithm leads to a generalization of Gessel's theorem to generalized Latin rectangles.

{\bf Input}

$\bullet$ Positive integer $k$, and another positive integer $N$ .

$\bullet$ $k(k-1)/2$ sets of {\it integers} $S_{ii'}$ ($1 \leq i < i' \leq k$) .

\vfill\eject

{\bf Output}

The first $N$ terms of the sequence

Number of $k \times n$ arrays $(M_{ij})$ where each row is a permutation of $\{1,2, \dots, n\}$, the first row is $1 \dots n$ and for all $1 \leq i < i' \leq k$, and  $1 \leq j < j' \leq n$, 

$$
M_{ij} \neq M_{i'j'} \quad, \quad if  \quad j'-j \in S_{ii'} \quad.
$$

The classical case of counting (reduced) $k \times n$ Latin rectangles corresponds to the case where all the $k(k-1)/2$ sets, $S_{ii'}$, happen to be the singleton set $\{0\}$.

The case $k=2$ , $ S_{12}=\{0,1\} $ is the famous (linear) {\it probl\`eme des m\'enages}, due to Eduard Lucas, and solved brilliantly, among others,
by  the human Irving Kaplansky [K].

One of the points we are trying to make is that using computers, with only a little bit of extra effort, one can treat a very general scenario.

{\bf A reminder about the Principle of Inclusion-Exclusion (our way)}

Suppose that you have a large set of {\it guys}, that is easy to count, and there is
a smaller set of {\it good guys} that is hard to count. In other words, the sum
$$
A:=\sum_{ g \in AllGuys} 1 \quad  ,
$$
is easy, but
$$
G:=\sum_{ g \in GoodGuys} 1 \quad  ,
$$
is hard.

We use two {\bf very deep} identities:

$\bullet$ $1+(-1)=0$.

$\bullet$ $0^i=0 \,\,\, if \,\,\, i>0$ \quad \quad , \quad  \quad $0^0=1 \quad .$

If some guy, let's call him Mr. $g$,  is {\bf not} a good guy, it means that he is a {\bf bad guy}, which means that
his set of {\bf sins}, $Sins(g)$,  is {\bf non-empty}, so we have
$$
G=\sum_{g \in AllGuys} 0^{|Sins(g)|} \quad  ,
$$
(the number of elements of a set $S$ is denoted by $|S|$). Now we use $0=1+(-1)$, and get
$$
G=\sum_{g \in AllGuys} (1+(-1))^{|Sins(g)|} \quad  .
$$
We need yet another deep identity. For any set $S$, we have:
\footnote{$^1$}
{\eightrm  \raggedright
Let $\chi(statement)$ equal $1$ or $0$ according to whether the statement is true or false. We have
$$
(1+(-1))^{|S|}= \prod_{s \in S} (1+(-1)) = \sum_{T \subset S} \prod_{s \in S} 1^{\chi (s \not \in T)}  (-1)^{\chi (s \in T)} 
$$
$$
=\sum_{T \subset S} \prod_{s \in S} (-1)^{\chi (s \in T)} =\sum_{T \subset S} \prod_{t \in T} (-1) 
=\sum_{T \subset S} (-1)^{|T|} \quad .
$$
}
$$
(1+(-1))^{|S|}= \sum_{T \subset S} (-1)^{|T|} \quad .
$$
Going back to the formula for $G$ we have
$$
G\, = \,
\sum_{g \in AllGuys} 0^{|Sins(g)|} =\sum_{g \in AllGuys} (1-1)^{|Sins(g)|} 
=\sum_{g \in AllGuys} \sum_{T \subset Sins(g)} (-1)^{|T|} \quad .
$$

So instead of straight-counting  the set of good guys, that is very hard to do, we do {\it weighted-counting} of the much larger set of {\bf pairs} $\{ [AnyGuy,S] \}$, were $S$ is a subset of the sins committed by $AnyGuy$,
and the {\it weight} is $(-1)^{|S|}$.

{\bf A quick reminder about Counting Derangements using the Principle of Inclusion-Exclusion} 

Before describing the  algorithm, let's revisit this old chestnut. We hope that our readers will not get offended, since
this will help  motivate the general case.

Instead of counting {\it good guys} (derangements) we weight-count the pairs $(\pi,S)$ where $\pi$ is {\it any} permutation of $\{1,2, \dots, n\}$ and $S$ is a subset of its {\it sins}, i.e.
a subset (possibly empty, possibly the whole thing) of the set of fixed points
$$
\{ 1 \leq i \leq n \, : \, \pi[i]=i \quad\} \quad,
$$
and give it weight $(-1)^{|S|}$. For example here is one such creature, where we indicate in {\bf boldface} the members of $S$
$$
\matrix{ 2&4&{\bf 3}&1&7&{\bf 6}&7&5&9\cr
         1&2&{\bf 3}&4&5&{\bf 6}&7&8&9} \quad,
$$
Here $S=\{3,6\}$ so its weight is $(-1)^2=1$.

Starting with such an  $S$, if it has $n-k$ elements,  and there are ${{n} \choose {n-k}}={{n} \choose {k}}$ such choices, each of weight $(-1)^{n-k}$, and the set of permutations $\pi$ that
can accompany any such $S$ {\bf definitely} has fixed points at the $n-k$ places, but all the other entries are at our disposal, provided that they
do not conflict with each other and the $n-k$ fixed members. There are $k!$ such $\pi$ that can serve as $S$-mates, and we get the famous formula for the number of derangements
$$
D_n \, = \, \sum_{k=0}^{n} \, (-1)^{n-k} \, {{n} \choose {k}}  k! \quad .
$$

Going back to the two-line notation, we can view the choice of $S$ as {\it tiling} of the  $2 \times n$   rectangle, $[0,1] \times [0,n-1]$, with horizontal translates of the three `tiles'
$$
\{[0,0]\} \quad, \quad \{[0,1]\} \quad, \quad \{[0,0],[0,1]\} \quad .
$$
So essentially what we did here is a two-step process.

{\bf 1.} Find the weight-enumerator of the set of tilings of $[0,1] \times [0,n-1]$ , where the weight of a tilings is the product of the formal variables  $t[s]$, where $s$ is a tile. In this simple case
the variables are
$$
t[\{[0,0]\}] \quad, \quad t[\{[0,1]\}] \quad, \quad t[\{[0,0],[0,1]\}] \quad .
$$

Since every column either has the two singleton tiles $\{[0,0]\}$, $\{[1,0]\}$ or the tile  $\{[0,0],[1,0]\}$, the weight-enumerator is
$$
\left ( \, t[\{[0,0]\}]\,t[\{[0,1]\}] \, + \, t[\{[0,0],[0,1]\}] \, \right )^n \quad .
$$

Given a tiling, how many ways can we put labels legally, so that the second row still has all distinct entries and the labels within each tile are the same?

Since the first row is fixed, Whenever there is a tile  $\{[0,0],[0,1]\}$, the corresponding entry in the second row is determined (it is the same as the one below it). As for
the cells that are singleton cells, those at the first row are determined, but those of the row can do what they want (as long as they are different from each other). So let's
replace $t[\{[0,0]\}]$ by $1$,  $t[\{[0,0], [0,1]\}]$ by $-1$ and  $t[\{[0,1]\}]$ by $x$, we get that the weight-enumerator according to the modified weight, is
$$
(x-1)^n \quad,
$$
and the final answer is obtained by applying the `umbra' $x^k \rightarrow k!$, that can be gotten by  $x^k \rightarrow \int_{0}^{\infty} x^k\, e^{-x} \, dx$,
that  leads to the alternative famous formula (probably due to Euler)
$$
D_n \, = \, \int_0^{\infty} (x-1)^n \, e^{-x} \, dx \quad .
$$

Note that the generating function of this weight-enumerator of tilings is the {\bf rational function}
$$
\frac{1}{1-z(1-x)} \quad ,
$$

Hence the ordinary generating function for the number of derangements is
$$
\sum_{n=0}^{\infty} \, D_n \, z^n \, = \, \int_0^{\infty} \frac{ e^{-x}}{1-z(1-x)} \, dx \quad .
$$

The above discussion may seem a bit too complicated for such a simple problem, and indeed it is, but we wanted to illustrate the general {\it approach} with the simplest possible not-entirely-trivial example.

{\bf The two-rowed case: The generalized M\'enages Problem}

This problem can be handled using {\it Rook polynomials} (see [Z3]), but here we present another way, that we believe is more natural.

Our problem now is

{\bf Input}

$\bullet$ A positive integer $N$

$\bullet$ a set $S$ of  integers.

{\bf Output}

The first $N$ terms of the sequence

``Number of permutation such that $i-\pi[i] \not \in S $''.

{\bf Algorithm}

{\bf 1.} For each $s \in S$ create a tile 
$$
\{ [0,0], [s,1] \} \quad,
$$
and also include the singleton tiles $\{[0,0]\}, \{[0,1]\} $.

{\bf 2.} Using {\it Symbolic Dynamical Programming} (as in [Z2]), adapted to the present case of {\it horizontal tiles}, implemented in the Maple package

{\tt https://sites.math.rutgers.edu/\~{}zeilberg/tokhniot/HorizTilings.txt} \quad,

we find the generating function of weight-enumerators for the boards $[0,1] \times [0,n-1]$ of the above tiles. Again for each tile $s$ with two members replace $t[s]$ by $-1$, and also
replace $t[\{[0,0]\}]$ by $1$, and   $t[\{[0,1]\}]$ by $x$. Finally, apply the umbra $x^k \rightarrow k!$. Equivalently, once Maple (automatically!) found
the rational function, $R(x,X)$, taylor-expand in powers of $X$, and then apply $x^k \rightarrow k!$, or equivalently find the first $N$ terms of
$$
\int_0^{\infty} \,e^{-x} \, R(x,X) \, dx \quad.
$$

We call $R(x,X)$ the {\it kernel}. This is implemented by procedure {\tt  Kernel2(S,x,X);}  in the Maple package {\tt GenLatinRecs.txt}, that can be downloaded from

{\tt https://sites.math.rutgers.edu/\~{}zeilberg/tokhniot/GenLatinRecs.txt} \quad,

where the set $S$ is as above. For example, typing

{\tt Kernel2($\{$0,1,-2$\}$,x,X);} gives

$$
\frac{\left(X +1\right) \left(X^{4}+\left(x -1\right) X^{2}+X +1 \right)}{X^{8}+\left(-x +1\right) X^{7}+\left(-2 x +2\right) X^{5}+\left(x^{2}-4 x +2\right) X^{4}+2 X^{2}+\left(-x +3\right) X +1} \quad .
$$

To get the first $N$ terms of the  sequence enumerating permutations such that $i-\pi[i] \not \in S$

type: {\tt GenDerSeq(S,N)}. For example to get the first $100$ straight M\'enage numbers, \hfill\break
A271 ( {\tt https://oeis.org/A000271}) , type

{\tt GenDerSeq($\{$0,1$\}$,100); } \quad .

To get the first $100$ terms of the sequence enumerating permutations such that $|\pi[i]-i| > 3$, type

{\tt GenDerSeq($\{$-3,-2,-1,0,1,2,3$\}$,100); } \quad ,

getting them in  one second. This is OEIS sequence A75852 {\tt https://oeis.org/A075852}, that currently(Aug. 13, 2021) only has $23$ terms.

The output file

{\tt https://sites.math.rutgers.edu/\~{}zeilberg/tokhniot/oGenLatinRecs1.txt} \quad ,

uses this procedure to generate the first $100$ terms of enumerating permutations such that $i- \pi[i] \not \in S$  for all $127$ non-empty subsets of $\{-3,-2,-1,0,1,2,3\}$.

{\bf The Three-Rowed case}

Now there are three kinds of {\bf forbidden} events.

$\bullet$ For  each $s \in S_{12}$ the set of corresponding `bad events' is $\{M_{1,j}=M_{2,j+s}\}$ for all $j$ such that both $j$ and $j+s$ are in $[1,n]$. These correspond to the translates of the edge $\{[0,0],[s,1]\}$.

$\bullet$ For  each $s \in S_{13}$ the set of corresponding `bad events' is $\{M_{1,j}=M_{3,j+s}\}$ for all $j$ such that both $j$ and $j+s$ are in $[1,n]$. These correspond to the translates of the edge $\{[0,0],[s,2]\}$.

$\bullet$ For  each $s \in S_{23}$ the set of corresponding `bad events' is $\{M_{2,j}=M_{3,j+s}\}$ for all $j$ such that both $j$ and $j+s$ are in $[1,n]$. These correspond to the translates of the edge $\{[0,1],[s,2]\}$.

In addition we have the {\bf singleton} tiles $\{[0,0]\}$ , $\{[0,1]\}$ , $\{[0,2]\}$.

Since the tiles are defined up to horizontal translation, we make the convention that the vertex with the smallest $y$-coordinate has $x$-coordinate $0$.

Every edge corresponds to a `violation'. Altogether they form a subgraph of $[0,n-1] \times [0,2]$, and naturally fall into {\bf connected components}, that
define a {\bf finite} set of {\bf tiles}, that the computer can find {\it all by itself} (in our Maple package {\tt GenLatinRecs.txt}). Then
using the method of [Z2], adapted to horizontal tilings (implemented in {\tt HorizTilings.txt} that is also included in the former package),
it {\it dynamically} sets up the finite-state system, and automatically sets up the corresponding set of linear equations for the weight-enumerators,
and then proceeds to solve them, {\it without a human touch}.

Now there are four kinds of tiles. Note that a tile can have at most one vertex from each of the three levels $y=0$, $y=1$, $y=2$, since all  rows are permutations,
otherwise we would have equalities within the same row.

$\bullet$ The singleton tiles $\{[0,0]\}$ , $\{[0,1]\}$ , $\{[0,2]\}$, with temporary weights  $t[\{[0,0]\}]$ , $t[\{[0,1]\}]$ , $t[\{[0,2]\}]$, to be adjusted later.

$\bullet$ Tiles that contain two vertices and one edge. There are three possibilities, those with $y$ coordinates $\{0,1\}$, those with $y$ coordinates $\{0,2\}$
those with $y$ coordinates in $\{1,2\}$, the weight of such a tile is $-t[s]$, (the $-1$ comes from the inclusion-exclusion for {\bf one} `sin').

$\bullet$ Tiles that contain three vertices and two edges, these have members from the three levels $y=0$, $y=1$, $y=2$.
The weight of such a tile is $t[s]$, (the sign, $(-1)^2$, comes from the inclusion-exclusion corresponding to the intersection of two sets).

$\bullet$ Tiles that contain three vertices and three edges, these also have members from the three levels $y=0$, $y=1$, $y=2$.
The weight of such a tile is $2\,t[s]$. The factor $2$ comes from the fact that now we have three underlying `sins'. Calling them $A,B,C$, we have
$$
|A \cap B| + |A \cap C| + |B \cap C| -|A \cap B \cap C| \quad ,
$$
giving the factor $(-1)^2+(-1)^2+(-1)^2+(-1)^3=2$.

At the end of the day, we now need four final variables for the {\it umbral finale}. Let's call them $x_1,x_2,x_3$ and $x_{23}$.
For all tiles $s$ that contain a vertex with $y=0$ replace $t[s]$ by $x_1$
(since the first row is assumed fixed, all their mates are fixed). Also replace $t[\{[0,0]\}]$ by $x_1$ (of course it is fixed).
The power of $x_1$ in a monomial will indicate the number of entries that are {\it already committed}.
For all the tiles that have a vertex with $y$-coordinate $1$ and $y$-coordinate $2$ (connecting the second and third rows) replace their $t[s]$ by the new formal variable $x_{23}$.
Finally replace $t[\{[0,1]\}]$ by $x_2$ and $t[\{[0,2]\}]$ by $x_3$.

Once the computer, {\it all by itself}, figured out, the set of participating tiles, and used {\it symbolic dynamical programming} to
find the generating function for the horizontal tilings with those weights, we need to perform the {\it umbral operation}:

$$
x_1^{a_1} \,  x_2^{a_2} \,  x_3^{a_3} \, x_{23}^{a_{23}} \rightarrow  {{n-a_1} \choose {a_{23}}} \cdot a_{23}! \cdot a_2! \cdot a_3! \quad .
$$

Let's explain it. There are $a_1$ labels totally committed from the tiles that contain a vertex with $y=0$. Then out of the remaining $n-a_1$ still-available labels 
we choose $a_{23}$  to label the tiles that connect the second and third row, and permute them in all possible ways. Now there are $a_2$ still-available labels for the second row and $a_3$ still-available labels for the third row, that
can be ordered as we please.

This is implemented in procedure {\tt GLR3seq(R12,R13,R23,N)}; where {\tt R12,R13,R23 } are the sets $S_{12},S_{13},S_{23}$, respectively, defined above. For seven examples see

{\tt https://sites.math.rutgers.edu/\~{}zeilberg/tokhniot/oGenLatinRecs2.txt} \quad.

A usual Latin $k \times n$ rectangle may be defined as an array of numbers from $\{1,2,\dots,n\}$ such that no two cells reachable via a Rook move are the same. Let's define a {\bf super Latin rectangle}
as one where no two cells reachable via a Queen move are the same.

For the $3 \times n$  super Latin rectangle, these conditions correspond to
$$
S_{12}=\{-1,0,1\} \quad,\quad
S_{23}=\{-1,0,1\} \quad,\quad
S_{13}=\{-2,0,2\} \quad.
$$

For the first 30 terms of this hard-to-compute sequence, enumerating $3 \times n$ {\bf super Latin rectangles}, see:

{\tt https://sites.math.rutgers.edu/\~{}zeilberg/tokhniot/oGenLatinRecs3.txt } \quad .

{\bf Generalizing Gessel's P-recursiveness theorem: For any fixed $k$ the number of Generalized $k \times n$ Latin rectangles is $P$-recursive in $n$}

Things get computationally very complicated for more than three rows, but the beauty of theory is that one can talk about things {\it in principle}.

The set of {\it violations} that can be committed by a $k \times n$ array, where every row is a permutation of $\{1,2,\dots n\}$  are
$$
M_{ij}=M_{i'j'} \quad, \quad 1 \leq i < i' \leq n \quad, \quad  j'-j \in S_{ii'} \quad .
$$
Each such violation corresponds to an {\bf edge} $\{[i,j],[i',j']\}$.  These edges naturally form {\bf connected components} that lead to a finite set of {\bf tiles}, that are defined up to horizontal translation.
Each such tile can contain at most $k$ vertices, at most one from each level $y=i$, $i=0, \dots, k-1$. 

For each such tile, $s$, look at all the subsets of the set of edges,(finitely many of them, of course) whose set of vertices consists of the set of vertices of that tile.
Now add $(-1)^{\#Edges}$, over all such subgraphs, let's call that number $c[s]$. This is the total contribution from the inclusion-exclusion that wind up at that very same tile.

After we get the generating function, in $X$, for the weight-enumerator of tilings according to the original variables $t[s]$, we replace
$t[s]$ by $c[s]\, t[s]$, to accommodate the fact that the context is inclusion-exclusion.

Now we switch to another set of variables. All the tiles that have a member of the first row, i.e. a vertex with $y$-coordinate $0$, are replaced by $x_1$. 
We also replace $t[\{[0,0]\}]$ by $x_1$.

Now for each subset $S$ of $\{1,2,\dots, k-1\}$ we introduce a new formal variable $x[S]$ and for each tile $s$  that has no $y$-coordinate at $y=0$, we replace $t[s]$ by $x[S]$ where $S$ is the
set of $y$-coordinates of the set of vertices of $s$.

Finally, in order, analogously to the three-rowed case, we do the `umbral substitution' and replace each monomial by the appropriate product of factorials and/or binomial coefficients.

It follows from the {\it general holonomic nonsense} of [Z1] that the resulting enumerating sequence is $P$-recursive.

{\bf Disclaimer}: The above is a not yet a full proof, but a sketch. Filling in the details would be rather daunting, and we prefer not to do it.

{\bf References}

[BL] K. P. Bogart and J. Q. Longyear, {\it Counting 3 by $n$  Latin rectangles},  Proc. Amer. Math. Soc. {\bf 54}(1976), 463-467.

[G1] Ira M. Gessel, {\it Counting three-line Latin rectangles}, Lect. Notes Math, {\bf 1234} (1986), Springer, 106-111.

[G2] Ira M. Gessel, {\it Counting Latin rectangles}, Bull. Amer. Math. Soc. (N.S.) {\bf 16} (1987), 79-82.

[K] Irving Kaplansky, {\it Solution of the probl\`eme des m\'enages},  Bull. Amer. Math. Soc. {\bf 49} (1943),  784-785. \hfill\break
[Reprinted in pp. 122-123, ``{\it Classical Papers in Combinatorics}'', edited by Ira Gessel and Gian-Carlo Rota, Birkhauser, 1987]

[R] John Riordan, `{\it Introduction to Combinatorial Analysis}', Dover, originally published by John Wiley, 1958.

[Z1] Doron Zeilberger, {\it A Holonomic systems approach to special functions
identities}, J. of Computational and Applied Math. {\bf 32}, 321-368 (1990). \hfill\break
{\tt  https://sites.math.rutgers.edu/\~{}zeilberg/mamarim/mamarimhtml/holonomic.html }\quad .

[Z2] Doron Zeilberger, {\it Automatic CounTilings},  The Personal Journal pf Shalosh B. Ekhad and Doron Zeilberger \hfill\break
{\tt https://sites.math.rutgers.edu/\~{}zeilberg/mamarim/mamarimhtml/tilings.html} \quad .

[Z3] Doron Zeilberger, {\it  Automatic Enumeration of Generalized M\'enage Numbers}, S\'eminaire Lotharingien de Combinatoire, {\bf B71a} (2014) \hfill\break
{\tt https://sites.math.rutgers.edu/\~{}zeilberg/mamarim/mamarimhtml/menages.html} \quad .

\bigskip
\hrule
\bigskip
George Spahn and Doron Zeilberger, Department of Mathematics, Rutgers University (New Brunswick), Hill Center-Busch Campus, 110 Frelinghuysen
Rd., Piscataway, NJ 08854-8019, USA. \hfill\break
Email: {\tt  gs828 at math dot rutgers dot edu} \quad, \quad {\tt DoronZeil] at gmail dot com}   \quad .

First Written: Fibonacci Day ({\bf Aug. 13, 2021}).

\end